\documentclass[12pt]{article}
\usepackage{graphicx}
\begin{document}

\newcommand{\Dh}{\hbox{\bf D}}
\newcommand{\Xh}{\hbox{\bf X}}
\newcommand{\Fh}{\hbox{\bf F}}


\title{\vskip.5cm
A limit-cycle solver for nonautonomous dynamical systems
\vskip.5cm}
\author{Rafael G. Campos$^*$ and Gilberto O. Arciniega\\ 
Facultad de Ciencias\\
Universidad de Colima\\
28045, Colima, Col., M\'exico.\\
\hbox{\small\tt rcampos@zeus.umich.mx, gilberto@cgic.ucol.mx}\\
}
\date{}
\maketitle
{\vskip.3cm
\noindent MSC: 65D25, 37M99, 37N05
}\\
\vspace*{1.5truecm}
\begin{center} Abstract \end{center}
A numerical technique used to solve boundary value problems is modified 
to find periodic steady-state solutions of nonautonomous dynamical 
systems. The technique uses a matrix representation of the time derivative
obtained through trigonometric interpolation of a periodic function. Such a 
differentiation matrix yields exact values for the derivative of a trigonometric 
polynomial and therefore, can be used as the main ingredient of a 
numerical method to solve nonlinear dynamical systems.
We apply this technique to obtain some limit cycles and bifurcation 
points of a sinusoidally driven pendulum and the steady-state response
of an electric circuit.
\vskip1.5cm
$^*$On leave from: Escuela de Ciencias F\'{\i}sico-Matem\'aticas, 
Universidad Michoacana, 58060 Morelia, Michoac\'an, M\'exico\\

\section{Introduction}
Dynamical systems can be given by a set of ordinary differential 
equations and some initial conditions. 
Such problems are usually solved by conventional 
procedures such as Runge-Kutta methods, multistep methods or
general linear methods (see for example \cite{Stu98}).
For instance, shooting methods \cite{Apr72},
extrapolation methods \cite{Ske80} and Newton methods based on the 
Poincare map \cite{Sem95} has been applied to accelerate the convergence
of the solution to the periodic steady-state of electric systems.
Since some of these implementations requires a considerable computational
effort, an alternative and more simple technique may be desired.\\
In a sequel of papers (see \cite{Cam00a} and references therein),
a Galerkin-collocation-type method for solving differential boundary value 
problems has been developed. In some restricted cases this technique yields 
exact results, i.e., the numerical output can be interpolated to obtain 
the functions that solve exactly the problem, and therefore it have been used 
to implement a numerical scheme to solve dynamical systems that
consists basically in the substitution of the derivatives
appearing in the differential equations by finite-dimensional matrices 
whose entries depend on certain set of points in a simple form. In some cases,
the nodes have to be chosen according to some criterion in order to
accelerate the convergence to the solution. The differentiation matrices used 
in this technique arise naturally in the context of the interpolation of 
functions. 
\\
The aim of this paper is to present a novel method to obtain
numerical approximants to the steady-state solutions of nonautonomous 
systems, based on the use of the differentiation matrix for periodic 
functions given in \cite{Cam00b}.\\
The main idea of the method is established in Sec. \ref{secdis} and
two examples are given in Sec. \ref{secap}. \\

\section{Discrete formalism} \label{secdis}
In this section we only present the main results of our discretization 
scheme for periodic functions; proofs and further applications can be 
found in \cite{Cam00b}.
Let $D$ be the $N\times N$ matrix whose components are given by
\[
D_{jk}=\cases{\displaystyle\mathop{{\sum}'}\limits_{l=1}^N
{1\over2}\cot{{t_j-t_l}\over 2},&{$j=k$},\cr\noalign{\vskip .5truecm}
\displaystyle {\tau_j\over{2\tau_k}}\csc{{t_j-t_k}\over 2}, &{$j\not=k$},\cr}
\]
where $-\pi<t_1<t_2<\cdots<t_N\le\pi$ and 
\[
\tau_j={d\over{dt}}\Big[\prod_{l=1}^N\sin{{t-t_l}\over 2}\Big]_{t_j}.
\]
Let $x(t)$ be a trigonometric polynomial of degree at most $n$ and 
$x$ and $x'$ denote the $N\times 1$ vectors whose elements 
are $x(t_j)$ and $x'(t_j)$, i.e., the values of the polynomial and its derivative
at the nodes respectively. Then, the matrix $D$ applied to $x$ becomes equal
to $x'$ whenever $N=2n+1$. Since the derivative of a trigonometric polynomial 
is again a trigonometric polynomial of the same degree, we have that  
\begin{equation}
x^{(k)}=D^k x,\quad  k=0,1,2,\ldots,
\label{xpk}
\end{equation}
where $x^{(k)}$ is the vector whose elements are given by the $k$th 
derivative of $x(t)$ evaluated at the nodes.
\\
The differentiation matrix $D$ takes the simple form 
\begin{equation}
D_{jk}=\cases{0,&{$j=k$},\cr\noalign{\vskip .5truecm}
\displaystyle {(-1)^{j+k}\over{2\sin{\pi\over N}(j-k) }}, &{$j\not=k$},\cr}
\label{djk}
\end{equation}
if the nodes are chosen to be the $N$ equidistant points\footnote{The 
same result is obtained for the set of points $t_j=\pi(2j-N-1)/N$.}
\begin{equation}
t_j=-\pi+{{2\pi j}\over N}, \quad j=1,2,\cdots,N.
\label{tj}
\end{equation}
If the periodic function $x(t)$ is not a polynomial, a residual vector depending 
on $x(t)$, $k$ and $N$ must be added to the right-hand side of (\ref{xpk}). 
However, it is expected that the norm of such a vector approaches zero as 
the number of nodes is increased since $x(t)$ can be expanded in a Fourier 
series. Therefore, if $N$ is great enough, the residual vector can be ignored 
and the function $x^{(k)}(t)$ can be approximated by an interpolation of the 
elements of $D^k x$.\\
Let us consider now a nonautonomous system of $m$ components described by 
\begin{equation}
\dot{x}=f(x,\omega t),
\label{es}
\end{equation}
where $x$ is the vector of components $x_1(t), x_2(t), \ldots , x_m(t)$ and
$f$ is a nonlinear vector function of $m+1$ variables, periodic in $t$, with period
$T=2\pi/\omega$. A very important 
problem is to find the response of the system in the steady-state regime. 
Taking into account the periodicity of this response, we can use the differentiation 
matrix (\ref{djk}) to obtain approximants to the steady-state solution of (\ref{es}) 
according to the following scheme.\\
First of all, we need to change $\omega t\rightarrow t$ in (\ref{es}) in order to
change the period of $f$ to $2\pi$. Thus, (\ref{es}) becomes
\begin{equation}
\omega \dot{x}=f(x,t),
\label{esp}
\end{equation}
Let us take an odd number $N$ of points $t_j$ as given by (\ref{tj}) and evaluate
(\ref{esp}) at each node to form an equality between $Nm\times 1$ vectors in 
such a way that the first $N$ entries of the left-hand side are the components
of the vector $\dot{x}_1$, i.e., $\dot{x}_1(t_1), \dot{x}_1(t_2),\ldots , \dot{x}_1(t_N)$; 
the following $N$ entries are the components of $\dot{x}_2$, i.e.,
$\dot{x}_2(t_1), \dot{x}_2(t_2),\ldots , \dot{x}_2(t_N)$, and so on. Now we can
approximate the vectors blocks $\dot{x}_k$ by $Dx_k$ to obtain the discrete form
of (\ref{esp}), which can be written as
\begin{equation}
\omega \sum_{l=1}^N D_{jl}x_k(t_l)=f_k(x_1(t_j),x_2,(t_j)\ldots,x_m(t_j),t_j),\qquad j=1,2,\ldots ,N,
\label{descom}
\end{equation}
where $k=1,2,\ldots, m$, or in the more compact form
\begin{equation}
\omega {\Dh}{\Xh}={\Fh}, \label{des}
\end{equation}
where ${\Fh}$ denotes the $Nm\times 1$ vector whose elements are given by the 
right-hand side of (\ref{descom}) first running $j$ and then $k$, ${\Xh}$ is the vector
whose elements are given by $x_k(t_j)$ (the unknown solution) ordered in a similar
way and ${\Dh}=1_m\otimes D$, where $1_m$ is the identity matrix of dimension $m$.
\\
The points on which our method is based are the following:
\begin{enumerate}
\item If (\ref{es}) has a limit cycle, then the solution of (\ref{des}) 
is an approximation of the steady-state solution of (\ref{es}). 
\item Since (\ref{des}) [or equivalently (\ref{descom})] is a system of $Nm$ 
nonlinear equations with $Nm$ unknowns $x_k(t_j)$, its solution can be
obtained by using a standard procedure (Newton's method for instance).
\end{enumerate}
This shows that a nonautonomous system in the steady-state regime can be described
approximately by a system of nonlinear algebraic equations. The procedure 
sketched in these statements is not concerned at all with the initial conditions of 
the system and yields simultaneously all the values of $x_k(t_j)$ (at all times). 
Thus, this method is quite different in essence to those designed as 
initial-value-problem solvers.\\
To obtain the solution of (\ref{descom}) we can use the Newton method or some 
variation of it. 
However, as it is well-known, not always is easy to give a good initial approximation
to attain convergence.  A good issue is to make a few Runge-Kutta integrations
of the system to yield an initial approximation ${\Xh}_0$. Since this procedure is
time-consuming, a more simple way to find ${\Xh}_0$ should be tried, but it will
depend on the problem to solve.
\\
Due to the nature of our method to approximate the steady-state solution of
(\ref{esp}), an algorithm for this technique will consist necessarily in the algorithm 
selected to solve the set of nonlinear equations. 
\\
In the following section we put this method at work on two important nonautonomous
dynamical systems.\\
\vskip.5cm


\section{Test cases} \label{secap}
We have chosen a chaotic mechanical system and a electric circuit of actual 
use as test cases.  To find their steady-state solutions, we rewrite the equations 
describing the dynamics of these problems in the discretized form (\ref{des}) 
and use a FORTRAN90 program and standard libraries running on a personal 
computer to solve them.\\
\subsection{Test case 1: the driven pendulum}
In spite that the classical pendulum is a very old problem, interest on it
is still growing. The pendulum becomes a chaotic system when it is driven at the
pivot point. Let us consider a rigid and planar pendulum consisting in a mass 
attached to a light rod of length $l$ which is vertically driven by a sinusoidal force
of the form $-A\cos\omega t$ and damped by a linear viscous force with damping
$\mu$ as in \cite{Bis99}. If $\theta$ denotes the 
angular displacement of the pendulum measured from the vertically downward 
position, the equation of motion is 
\begin{equation}
{{d^2\theta}\over{dt^2}}+a{{d\theta}\over{dt}}+(1+b\cos\omega t)\sin\theta=0,
\label{eqpen}
\end{equation}
where
\[
a={{2\mu}\over{\sqrt{lg}}},\qquad b={{A\omega^2}\over{l}}.
\]
To compare the steady-state solutions yielded by our procedure with those 
obtained by other authors we take $0\le b\le 200$ 
(the driving amplitude) as the control parameter, $a=0.1$ and $\omega=17.5$.
For values of $b$ near zero, a good initial guess is a vector ${\Xh}_0$ whose entries 
are close to $\pi$. To catch more solutions for a given value of $b$, one can try
a vector with a sinusoidal deviation from $\pi$. Once we have found a solution ${\Xh}(p)$, 
we proceed to obtain ${\Xh}(p+\delta p)$ by using ${\Xh}_0={\Xh}(p)$. We have taken 
$N=101$ nodes of the form (\ref{tj}) in this case.
\\
As it is known, this system presents period-doubling and bifurcation. Figure 1 displays
these phenomena for solutions of (\ref{eqpen}) close to the inverted state $\theta=\pi$.
To make the drawing more simple, we only plot the maximum and minimum values of 
the angular displacement $\theta(t)$ as functions of the control parameter $b$. Other 
solutions, as the hanging state, are not considered here.
\\ 
The displayed information is in line with the findings of other authors 
(see for example \cite{Bis99}), but in our case the problem of stability of the solution is 
not revised. In Fig. 1, $\theta_1$ is the solution corresponding to the inverted state,
$\theta_2$ and $\theta_3$ correspond in essence to the same period-1 solution (one can
get the second from the first through a reflection on the axis $\theta=\pi$) and
$\theta_4$ and $\theta_5$ are period-2 solutions oscillating both of them about the vertical.
The  asymmetrical way in which $\theta_5$ oscillates is illustrated in Fig. 2.\\
\vskip1cm
\hbox to \textwidth{\hfill\scalebox{0.6}{\includegraphics{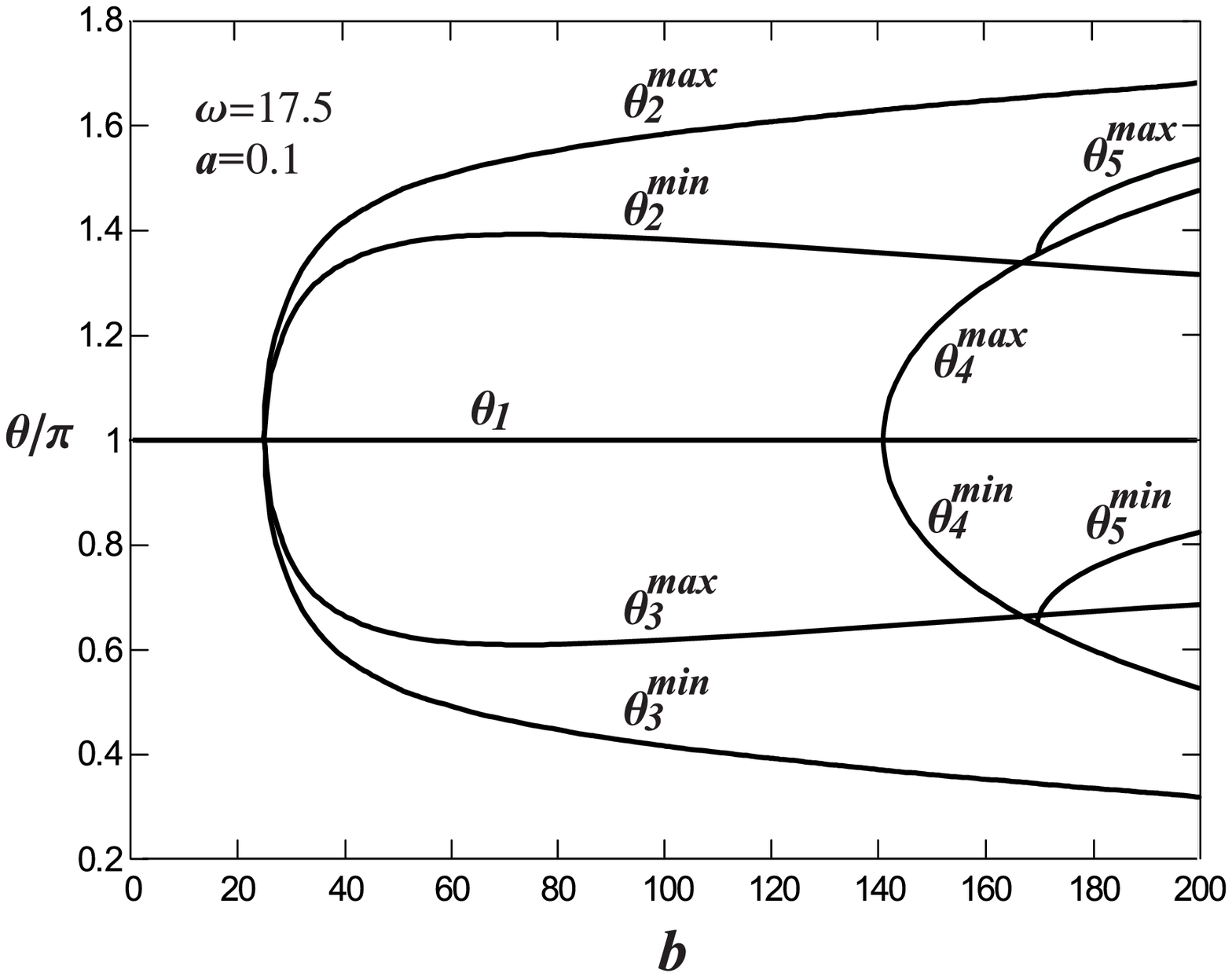}}\hfill}
\begin{center}
\begin{minipage}{14cm}
{\small 
Figure 1: Maximum and minimum values of the angular displacement $\theta$ of 
the vertically driven pendulum (\ref{eqpen}) vs. the driven amplitude $b$ for $a=0.1$ 
and $\omega=17.5$. For any value of $b$, the solution $\theta_1=\pi$ (the inverted
pendulum) always was found. The curves labeled by $\theta_k$, $k=2,\cdots,5$
correspond to different solutions.}
\end{minipage}
\end{center}
\vskip1cm
\hbox to \textwidth{\hfill\scalebox{0.6}{\includegraphics{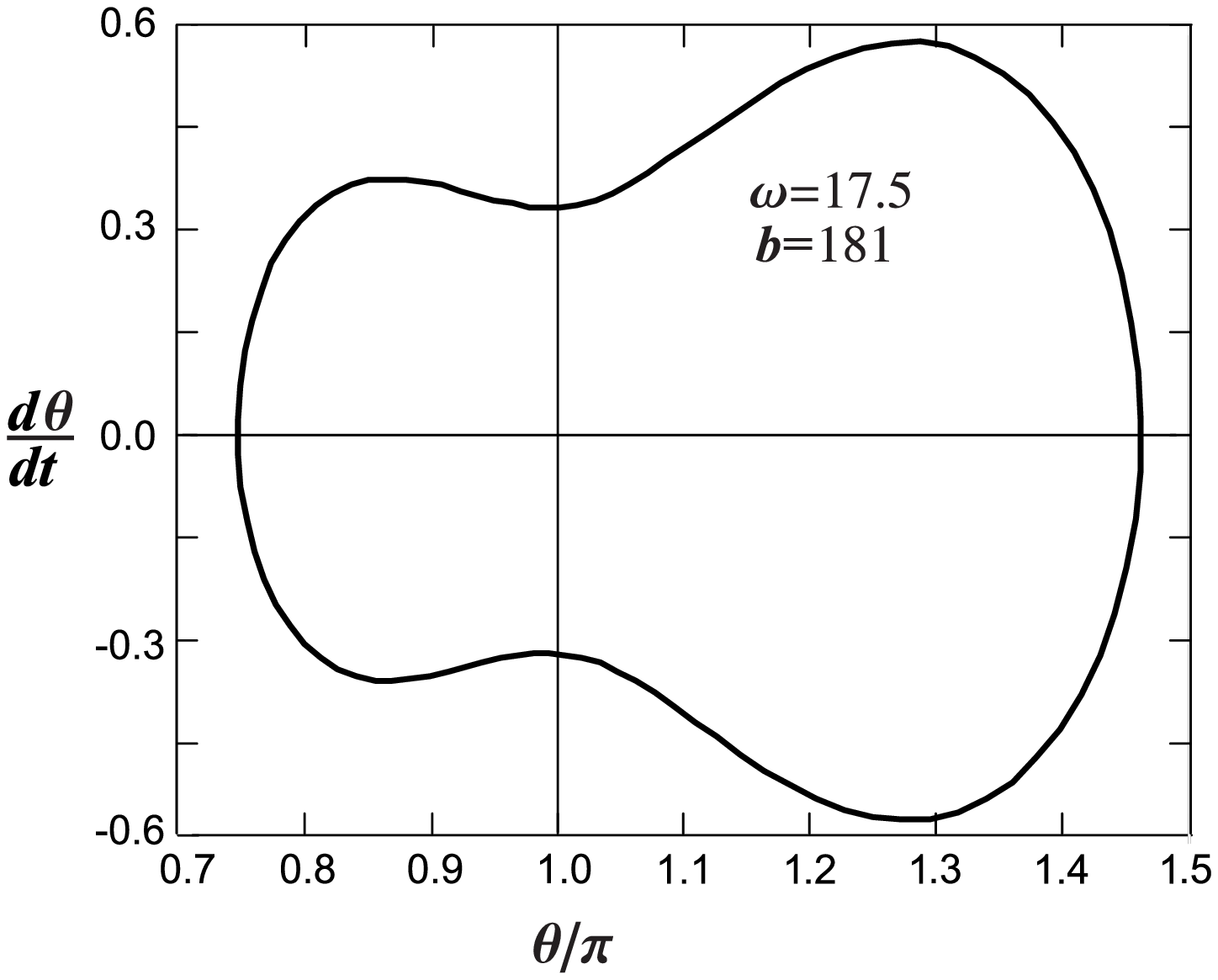}}\hfill}
\begin{center}
\begin{minipage}{14cm}
{\small 
Figure 2: Phase-space diagram of the period-2 "mariachi" solution $\theta_4$ 
for $a=0.1$, $\omega=17.5$ and $b=181$. }
\end{minipage}
\end{center}
\vskip1cm

\subsection{Test case 2: a commutation circuit}
The actual interest in circuits as the one depicted in Fig. 3, relies on the extended
use of them in modern power supplies. Usually a primary DC power supply 
(a battery or the rectified and filtered 60 Hz main line) is
commutated at high frequencies (tens of KHz to MHz) in order to obtain a pulse 
waveform, which can be coupled to the load in some way to obtain a
predetermined load voltage. On the circuit shown in Fig. 3, $V_s$ is a simplified 
model of a transformer secondary, which is followed by a rectifier-filter circuit and a 
resistive load ($R_4$). $R_1$ stands 
for the series resistance of the diode and $R_2$ and $R_3$ model the equivalent series 
resistance of $C_1$ and $C_2$ respectively. $R_2$ and $R_3$ are important
parameters because they increase ripple voltages and circuit losses.\\

\vskip1cm
\hbox to \textwidth{\hfill\scalebox{0.6}{\includegraphics{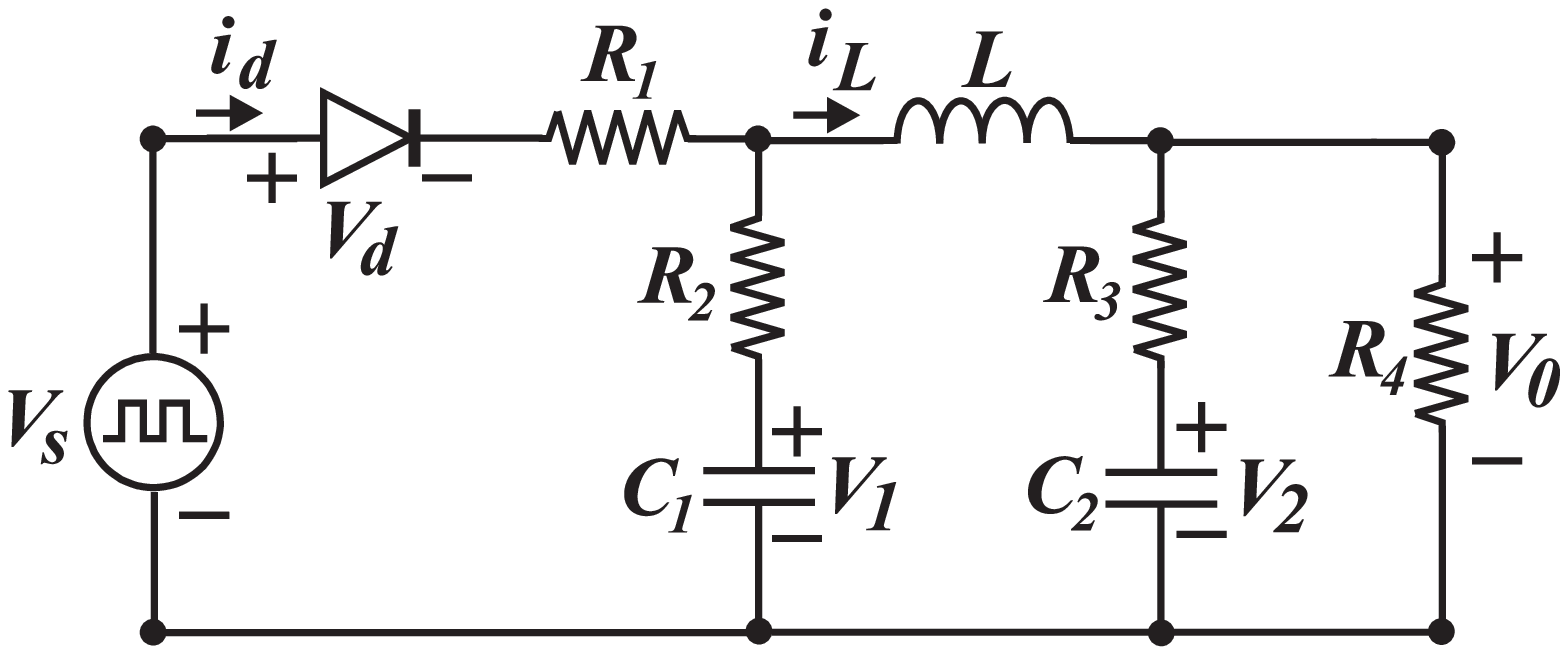}}\hfill}
\begin{center}
\begin{minipage}{14cm}
{\small 
\begin{center} Figure 3: A typical commutation circuit. \end{center} }
\end{minipage}
\end{center}
\vskip1cm
The equations for the state variables $V_1=x_1$, $V_2=x_2$ and $i_{L}=x_3$ 
can be simplified to the form 
\begin{eqnarray}
\dot{x_1}&=&{1\over{C_1R_2}}\Big[V_s-x_1-V_d-i_sR_1\Big(\exp({{V_dq}\over{\eta kT}})-1\Big)\Big], \nonumber\\
\dot{x_2}&=&{1\over{C_2(R_3+R_4)}}\Big(-x_2+R_4x_3\Big), \label{eqcir}\\
\dot{x_3}&=&{1\over{L}}\Big[V_s-{R_4\over{(R_3+R_4)}}x_2-{R_3R_4\over{(R_3+R_4 )}}x_3
-V_d-i_sR_1\Big(\exp({{V_dq}\over{\eta kT}})-1\Big)\Big]. \nonumber
\end{eqnarray}
where 
\[
V_d=V_s-C_1(R_1+R_2)\dot{x_1}-x_1-R_1x_3,
\]
$k$ is the Boltzmann constant, $q$ is the electron charge, $T$ the absolute temperature 
which is taken at the standard value $T=300^\circ K$ and $\eta$ is the emission
coefficient. The driving voltage is chosen 
as a square wave pulse, i.e., $V_s(t)=A_m\hbox{sgn}(2\pi t/T)$, $t\in [-T/2,T/2]$,
and the remaining parameters of this problem are taken at typical values: $A_m=5.6\hbox{V}$, 
$T=10^{-5}\hbox{s}$, $i_s=10^{-8}\hbox{A}$, $R_1=0.0149\Omega$, $R_2=0.15\Omega$, 
$R_3=0.2\Omega$, $R_4=2.0\Omega$, 
$C_1=470.0\times10^{-6}\hbox{F}$, $C_2=20.0\times10^{-6}\hbox{F}$, 
$L=20.0\times10^{-6}\hbox{H}$, $\eta=0.8953$.\\
In order to have results to compare with, we simulated the circuit using a LT version 
of the Spice program which is currently available for free from Linear Technology, 
a major semiconductor devices manufacturer. Spice \cite{Nag73} is a standard tool used for 
electric circuits simulation. Such a program compute the steady-state response in 
a sequential way and therefore, through a number of periods of the excitation.
For this circuit, the response of the system was computed by the Spice simulation 
over 150 cycles to achieve the steady-state. A maximum step size of $0.04 \mu\hbox{s}$
was specified. In our case, we use 251 equispaced time values $t_j$ of the form (\ref{tj}) 
(yielding a constant step size approximately equal to $0.04 \mu\hbox{s}$)
and all of the steady-state response values were computed at once. 
Since Spice uses an adaptive strategy for the step size, it is not possible to compare 
both of the results in a point-by-point basis. However, it comes out that they agree 
globally up to $10^{-2}$. \\
The output yielded by the present method and the Spice simulation for the more
used variables $i_d=x_3+C_1\dot{x_1}$ and $V_0=C_2R_3\dot{x_2}+x_2$,
is displayed in Figs. (4a) and (4b), respectively. 
\vskip1cm
\hbox to \textwidth{\hfill\scalebox{0.8}{\includegraphics{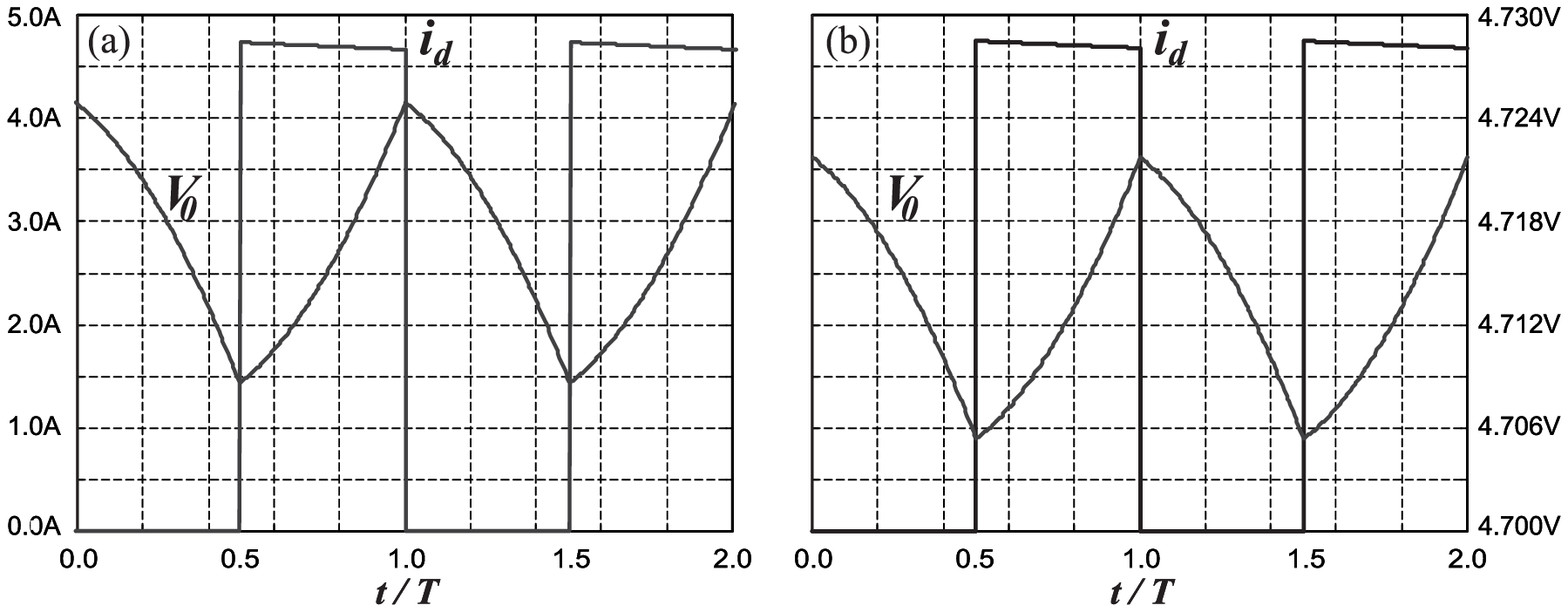}}\hfill}
\begin{center}
\begin{minipage}{14cm}
{\small 
Figure 4: Steady-state solution for the diode current  $i_d=x_3+C_1\dot{x_1}$
and the voltage $V_0=C_2R_3\dot{x_2}+x_2$ along two periods obtained by
the present method with $N=251$ nodes (4a) and by a Spice simulation (4b). 
The left-hand side and right-hand side vertical axes are used for different scales.}
\end{minipage}
\end{center}


\section{Final remarks} \label{conclus}
As the simple tests of the previous section show, the method presented in this
paper gives a novel approach to obtain the limit cycles of nonautonomous 
systems and represents an alternative to the conventional methods of
integration. Based on a well-behaved differentiation matrix for periodic functions,
it departs from the standard procedures in the way in which the limit cycle is found:
the usual integration of the system for long times is replaced by the problem of
finding the solution of a set of nonlinear algebraic equations. This feature makes
this method a tool for studying dynamical systems from a new point of view.
\\


\section{Acknowledgments} 
RGC wants to thank Dr. A. Medina and N. Garc\'{\i}a for very useful discussions
and suggestions.

\end{document}